\documentclass[12pt]{article}

\usepackage{amsmath,amssymb,amsthm,amscd,a4wide,upref}
\usepackage[toc]{appendix}

\usepackage{hyperref}
\hypersetup{colorlinks,citecolor=blue,filecolor=black,linkcolor=blue,urlcolor=blue}
\usepackage{enumerate}
\usepackage{mathtools}

\begin{document}

\newcommand{\ad}{{\rm ad}}
\newcommand{\cri}{{\rm cri}}
\newcommand{\End}{{\rm{End}\ts}}
\newcommand{\Rep}{{\rm{Rep}\ts}}
\newcommand{\res}{{\rm{res}}}
\newcommand{\Hom}{{\rm{Hom}}}
\newcommand{\Mat}{{\rm{Mat}}}
\newcommand{\ch}{{\rm{ch}\ts}}
\newcommand{\chara}{{\rm{char}\ts}}
\newcommand{\diag}{{\rm diag}}
\newcommand{\non}{\nonumber}
\newcommand{\wt}{\widetilde}
\newcommand{\wh}{\widehat}
\newcommand{\ot}{\otimes}
\newcommand{\la}{\lambda}
\newcommand{\La}{\Lambda}
\newcommand{\De}{\Delta}
\newcommand{\al}{\alpha}
\newcommand{\be}{\beta}
\newcommand{\ga}{\gamma}
\newcommand{\Ga}{\Gamma}
\newcommand{\ep}{\epsilon}
\newcommand{\ka}{\kappa}
\newcommand{\vk}{\varkappa}
\newcommand{\vt}{\vartheta}
\newcommand{\si}{\sigma}
\newcommand{\vs}{\varsigma}
\newcommand{\vp}{\varphi}
\newcommand{\de}{\delta}
\newcommand{\ze}{\zeta}
\newcommand{\om}{\omega}
\newcommand{\Om}{\Omega}
\newcommand{\ee}{\epsilon^{}}
\newcommand{\su}{s^{}}
\newcommand{\hra}{\hookrightarrow}
\newcommand{\ve}{\varepsilon}
\newcommand{\ts}{\,}
\newcommand{\pr}{^{\tss\prime}}
\newcommand{\vac}{\mathbf{1}}
\newcommand{\di}{\partial}
\newcommand{\qin}{q^{-1}}
\newcommand{\tss}{\hspace{1pt}}
\newcommand{\Sr}{ {\rm S}}
\newcommand{\U}{ {\rm U}}
\newcommand{\BL}{ {\overline L}}
\newcommand{\BE}{ {\overline E}}
\newcommand{\BP}{ {\overline P}}
\newcommand{\AAb}{\mathbb{A}\tss}
\newcommand{\CC}{\mathbb{C}\tss}
\newcommand{\KK}{\mathbb{K}\tss}
\newcommand{\QQ}{\mathbb{Q}\tss}
\newcommand{\SSb}{\mathbb{S}\tss}
\newcommand{\TT}{\mathbb{T}\tss}
\newcommand{\ZZ}{\mathbb{Z}\tss}
\newcommand{\DY}{ {\rm DY}}
\newcommand{\X}{ {\rm X}}
\newcommand{\Y}{ {\rm Y}}
\newcommand{\Z}{{\rm Z}}
\newcommand{\Ac}{\mathcal{A}}
\newcommand{\Lc}{\mathcal{L}}
\newcommand{\Mc}{\mathcal{M}}
\newcommand{\Pc}{\mathcal{P}}
\newcommand{\Qc}{\mathcal{Q}}
\newcommand{\Rc}{\mathcal{R}}
\newcommand{\Sc}{\mathcal{S}}
\newcommand{\Tc}{\mathcal{T}}
\newcommand{\Bc}{\mathcal{B}}
\newcommand{\Ec}{\mathcal{E}}
\newcommand{\Fc}{\mathcal{F}}
\newcommand{\Gc}{\mathcal{G}}
\newcommand{\Hc}{\mathcal{H}}
\newcommand{\Uc}{\mathcal{U}}
\newcommand{\Vc}{\mathcal{V}}
\newcommand{\Wc}{\mathcal{W}}
\newcommand{\Yc}{\mathcal{Y}}
\newcommand{\Ar}{{\rm A}}
\newcommand{\Br}{{\rm B}}
\newcommand{\Ir}{{\rm I}}
\newcommand{\Fr}{{\rm F}}
\newcommand{\Jr}{{\rm J}}
\newcommand{\Or}{{\rm O}}
\newcommand{\GL}{{\rm GL}}
\newcommand{\Spr}{{\rm Sp}}
\newcommand{\Rr}{{\rm R}}
\newcommand{\Zr}{{\rm Z}}
\newcommand{\ZX}{{\rm ZX}}
\newcommand{\gl}{\mathfrak{gl}}
\newcommand{\middd}{{\rm mid}}
\newcommand{\ev}{{\rm ev}}
\newcommand{\Pf}{{\rm Pf}}
\newcommand{\Norm}{{\rm Norm\tss}}
\newcommand{\oa}{\mathfrak{o}}
\newcommand{\spa}{\mathfrak{sp}}
\newcommand{\osp}{\mathfrak{osp}}
\newcommand{\f}{\mathfrak{f}}
\newcommand{\g}{\mathfrak{g}}
\newcommand{\h}{\mathfrak h}
\newcommand{\bgot}{\mathfrak b}
\newcommand{\n}{\mathfrak n}
\newcommand{\z}{\mathfrak{z}}
\newcommand{\Zgot}{\mathfrak{Z}}
\newcommand{\p}{\mathfrak{p}}
\newcommand{\sll}{\mathfrak{sl}}
\newcommand{\agot}{\mathfrak{a}}
\newcommand{\qdet}{ {\rm qdet}\ts}
\newcommand{\Ber}{ {\rm Ber}\ts}
\newcommand{\HC}{ {\mathcal HC}}
\newcommand{\cdet}{{\rm cdet}}
\newcommand{\rdet}{{\rm rdet}}
\newcommand{\tr}{ {\rm tr}}
\newcommand{\gr}{ {\rm gr}\ts}
\newcommand{\str}{ {\rm str}}
\newcommand{\loc}{{\rm loc}}
\newcommand{\Gr}{{\rm G}}
\newcommand{\sgn}{ {\rm sgn}\ts}
\newcommand{\sign}{{\rm sgn}}
\newcommand{\ba}{\bar{a}}
\newcommand{\bb}{\bar{b}}
\newcommand{\eb}{\bar{e}}
\newcommand{\bi}{\bar{\imath}}
\newcommand{\bj}{\bar{\jmath}}
\newcommand{\bk}{\bar{k}}
\newcommand{\bl}{\bar{l}}
\newcommand{\hb}{\mathbf{h}}
\newcommand{\Sym}{\mathfrak S}
\newcommand{\fand}{\quad\text{and}\quad}
\newcommand{\Fand}{\qquad\text{and}\qquad}
\newcommand{\For}{\qquad\text{or}\qquad}
\newcommand{\OR}{\qquad\text{or}\qquad}
\newcommand{\grpr}{{\rm gr}^{\tss\prime}\ts}
\newcommand{\degpr}{{\rm deg}^{\tss\prime}\tss}

\numberwithin{equation}{section}

\newtheorem{thm}{Theorem}[section]
\newtheorem{lem}[thm]{Lemma}
\newtheorem{prop}[thm]{Proposition}
\newtheorem{cor}[thm]{Corollary}
\newtheorem{conj}[thm]{Conjecture}
\newtheorem*{mthm}{Main Theorem}
\newtheorem*{mthma}{Theorem A}
\newtheorem*{mthmb}{Theorem B}
\newtheorem*{mthmc}{Theorem C}
\newtheorem*{mthmd}{Theorem D}

\theoremstyle{definition}
\newtheorem{defin}[thm]{Definition}

\theoremstyle{remark}
\newtheorem{remark}[thm]{Remark}
\newtheorem{example}[thm]{Example}

\newcommand{\bth}{\begin{thm}}
\renewcommand{\eth}{\end{thm}}
\newcommand{\bpr}{\begin{prop}}
\newcommand{\epr}{\end{prop}}
\newcommand{\ble}{\begin{lem}}
\newcommand{\ele}{\end{lem}}
\newcommand{\bco}{\begin{cor}}
\newcommand{\eco}{\end{cor}}
\newcommand{\bde}{\begin{defin}}
\newcommand{\ede}{\end{defin}}
\newcommand{\bex}{\begin{example}}
\newcommand{\eex}{\end{example}}
\newcommand{\bre}{\begin{remark}}
\newcommand{\ere}{\end{remark}}
\newcommand{\bcj}{\begin{conj}}
\newcommand{\ecj}{\end{conj}}

\newcommand{\bal}{\begin{aligned}}
\newcommand{\eal}{\end{aligned}}
\newcommand{\beq}{\begin{equation}}
\newcommand{\eeq}{\end{equation}}
\newcommand{\ben}{\begin{equation*}}
\newcommand{\een}{\end{equation*}}

\newcommand{\bpf}{\begin{proof}}
\newcommand{\epf}{\end{proof}}

\def\beql#1{\begin{equation}\label{#1}}

\title{\Large\bf Higher order Hamiltonians for
 the trigonometric Gaudin model}

\author{{Alexander Molev\quad and\quad Eric Ragoucy}}

\date{} 
\maketitle

\vspace{5 mm}

\begin{abstract}
We consider the trigonometric classical $r$-matrix for $\gl_N$
and the associated quantum Gaudin model. We produce higher Hamiltonians
in an explicit form
by applying the limit $q\to 1$ to elements of the Bethe subalgebra
for the $XXZ$ model.

\bigskip

Preprint LAPTH-004/18

\end{abstract}

\vspace{10 mm}

\noindent
School of Mathematics and Statistics\newline
University of Sydney,
NSW 2006, Australia\newline
alexander.molev@sydney.edu.au

\vspace{7 mm}

\noindent
Laboratoire de Physique Th\'{e}orique LAPTh,
CNRS and Universit\'{e} de Savoie\newline
BP 110, 74941 Annecy-le-Vieux Cedex, France\newline
eric.ragoucy@lapth.cnrs.fr

\bigskip

\section{Introduction}
\label{sec:int}

Explicit higher Hamiltonians for the rational Gaudin model
associated with $\gl_N$ were produced by Talalaev~\cite{t:qg} by making use
of the Bethe subalgebra of the Yangian for $\gl_N$ and taking
a classical limit; see also \cite{mtv:be}.
Some related families of higher Hamiltonians
and their analogues for the orthogonal and symplectic Lie algebras
were produced by using the {\em center at the critical level}
following the general approach of Feigin, Frenkel and Reshetikhin~\cite{ffr:gm};
see \cite{cf:mm}, \cite{cm:ho} and also \cite{m:so} for more details and references.
In particular, such a family arises from the
coefficients of the differential operators
$\tr\big(\di_u+E(u)\big)^k$ with $k=1,2,\dots$ which
form a commutative subalgebra of $\U\big(t^{-1}\gl_N[t^{-1}]\big)$. Here
$E(u)=[E_{ij}(u)]$ is the matrix with the entries
\ben
E_{ij}(u)=\sum_{n=0}^{\infty} E_{ij}[-n-1]\tss u^n,\qquad E_{ij}[-n-1]=E_{ij}\tss t^{-n-1},
\een
where we use the standard basis elements $E_{ij}$ of $\gl_N$.
The commuting family of the coefficients of the power series $\tr\ts E(u)^2$ can be regarded as
a generating series of the quadratic Gaudin Hamiltonians as considered by Sklyanin~\cite{s:sv}.

Both rational and trigonometric Gaudin models were studied by Jur\v{c}o~\cite{j:cy}.
They are associated with the
corresponding {\em classical $r$-matrix} $r(x)$. In the trigonometric case
$t^{-1}\gl_N[t^{-1}]$ is replaced by the extended Lie algebra
$\wh\g^{\tss +}=\bgot^+\oplus t^{-1}\gl_N[t^{-1}]$,
where $\bgot^+$ is the subalgebra of $\gl_N$ spanned by the elements $E_{ij}$ with $i\leqslant j$.
Accordingly, $E(u)$ is replaced by the matrix $\Lc^+(u)=[\Lc^+_{ij}(u)]$ with
\beql{lplus}
\Lc^+_{ij}(u)=\sum_{n=0}^{\infty} \Lc^+_{ij}[-n]\tss u^n,
\eeq
where $\Lc^+_{ij}[-n]=-2\tss E_{ij}\tss t^{-n}$ for $n\geqslant 1$ and
$\Lc^+_{ij}[0]=-\big(1+\sign(j-i)\big) E_{ij}$ assuming $\sign(0)=0$.
It is shown in \cite{j:cy} that the coefficients of the series $\tr\ts\Lc^+(u)^2$
are pairwise commuting elements of $\U(\wh\g^{\tss +})$. Moreover,
as with the rational Gaudin model, this series plays the role of the generating function
for quadratic Hamiltonians. Namely, taking the image in the tensor product
of the vector representations, we get
\beql{lrep}
\Lc^+(u)\mapsto r_{01}(u/a_1)+\dots+r_{0l}(u/a_l)
\eeq
for some parameters $a_i$, where
\beql{rx}
r(x)=\sum_{i,j=1}^N \Big(\ts\frac{1+x}{1-x}+\sign(j-i)\Big)\ts e_{ij}\ot e_{ji}
\eeq
is a {\em trigonometric classical $r$-matrix}. It satisfies the {\em classical
Yang--Baxter equation}
\ben
\big[r_{12}(x_1/x_2),r_{23}(x_2/x_3)\big]+\big[r_{23}(x_2/x_3),r_{31}(x_3/x_1)\big]+
\big[r_{31}(x_3/x_1),r_{12}(x_1/x_2)\big]=0
\een
together with the {\em skew-symmetry condition}
\ben
r_{12}(x)+r_{21}(1/x)=0.
\een
Taking
the residue at $a_i$, we recover the $i$-th Gaudin Hamiltonian
\beql{quadham}
\underset{u=a_i}{\res}\ts\tr\ts\Lc^+(u)^2=2\tss a_i\ts \sum_{j\ne i} r_{ij}(a_i/a_j),
\eeq
assuming the parameters $a_i$ are all distinct and nonzero.

Our main result is a construction of higher order Hamiltonians for the
trigonometric Gaudin model. They are obtained from a commuting family of elements of
$\U(\wh\g^{\tss +})$ which occur as the coefficients of formal series
written explicitly in terms of $\Lc^+(u)$. This commuting family is analogous to
the one produced from the differential operators $\tr\big(\di_u+E(u)\big)^k$
because the highest degree term of the corresponding operator
coincides with $\tr\ts\Lc^+(u)^k$.
By using the representation \eqref{lrep} one gets a commuting family of
higher order Hamiltonians which are pairwise commuting operators in the
tensor product of the vector representations.

\medskip

We are grateful to Nicolas Cramp\'{e} for many useful discussions.
This work was completed during the first author's visit to
the {\it Laboratoire d'Annecy-le-Vieux de Physique Th\'eorique\/}.
He thanks the lab for the support and hospitality.

\section{Trigonometric Gaudin model}
\label{sec:tgm}

The commutation relations of the Lie algebra $\wh\g^{\tss +}$ admit the matrix form
\ben
\big[\Lc^+_1(u),\Lc^+_2(v)\big]=\big[\Lc^+_1(u)+\Lc^+_2(v),r_{12}(u/v)\big],
\een
where both sides take values in $\End\CC^N\ot \End\CC^N\ot \U(\wh\g^{\tss +})$ and
the subscripts indicate the copies of the endomorphism algebra.
For all $s\geqslant 1$ consider the multiple tensor products
\beql{mtp}
\underbrace{\End\CC^N\ot \dots\ot\End\CC^N}_{s}\ot\tss \U(\wh\g^{\tss +}).
\eeq
Introduce the function $T(y)$ in a variable $y$ with values
in $\End\CC^N\ot \End\CC^N$ by
\beql{ty}
T(y)=\sum_{i=1}^N  e_{ii}\ot e_{ii}
+\frac{1}{1-y}\ts\sum_{i<j} e_{ij}\ot e_{ji}
+\frac{1}{1+y}\ts\sum_{i>j} e_{ij}\ot e_{ji}.
\eeq
For any $1\leqslant a<b\leqslant s$ we let $T_{ab}(y)$ denote the function $T(y)$
regarded as an element of \eqref{mtp} associated with the $a$-th and $b$-th copies
of $\End\CC^N$ and as the identity element in all the remaining tensor factors.
Now define differential operators $\theta_m\in \U(\wh\g^{\tss +})[[u,\di_u]]$ by means of the generating
function
\beql{thetham}
\sum_{m=1}^{\infty} \theta_m\tss y^m=\sum_{s=1}^{\infty}\tss y^s\ts
\tr^{}_{1,\dots,s}\ts T_{s-1\ts s}(y)\dots T_{1\ts 2}(y)\ts\Lc_1\dots \Lc_s,
\eeq
where $\Lc=2\tss u\tss \di_u-\Lc^+(u)$ and the trace is taken over all $s$ copies
of $\End\CC^N$. We can write
$T(y)$ as the series
\beql{tyexp}
T(y)=P+\sum_{r=1}^{\infty} \big(\tss\overline{\Tc}y^{2r}+\Tc y^{2r+1}\big),
\eeq
where
\beql{permut}
P=
\sum_{i,j=1}^N  e_{ij}\ot e_{ji},\qquad
\Tc=
\sum_{i,j=1}^N \sign(j-i)\ts e_{ij}\ot e_{ji}\Fand
\overline{\Tc}=
\sum_{i\ne j} e_{ij}\ot e_{ji}.
\eeq
Note that $\Tc=r(-1)$ is the value of the classical $r$-matrix \eqref{rx} at $x=-1$.
Taking the coefficient of $y$ in \eqref{thetham} we get
\ben
\theta_1=\tr\ts\Lc=2\tss N\tss u\tss \di_u-\tr\ts \Lc^+(u).
\een
The coefficients of the series $\tr\ts \Lc^+(u)$ are central in $\U(\wh\g^{\tss +})$.
Furthermore,
\begin{multline}
\theta_2=\tr_{1,2}\ts P_{12}\Lc_1\Lc_2=\tr\ts \Lc^2
=\tr \big(2\tss u\tss \di_u-\Lc^+(u)\big)^2\\[0.4em]
{}=4\tss N\tss u^2\tss \di_u^{\tss 2}
-4\tss u\big(\tr\ts\Lc^+(u)-N\big)\di_u-2\tss u\ts\tr\ts \Lc^+(u)'+\tr\ts \Lc^+(u)^2
\non
\end{multline}
and
\begin{multline}
\theta_3=\tr_{1,2,3}\ts P_{23}\tss P_{12}\Lc_1\Lc_2\Lc_3+\tr_{1,2}\ts \Tc_{12}\Lc_1\Lc_2
=\tr\ts \Lc^3+\tr_{1,2}\ts \Tc_{12}\Lc_1\Lc_2\\[0.4em]
{}=\tr \big(2\tss u\tss \di_u-\Lc^+(u)\big)^3+\sum_{i,j=1}^N \sign(i-j)\Lc^+_{ij}(u)\Lc^+_{ji}(u).
\non
\end{multline}
For any $m\geqslant 1$ the differential operator $\theta_m$ takes the form
\ben
\theta_m=\theta^{(0)}_m\tss \di_u^{\tss m}+\dots+\theta^{(m-1)}_m\tss\di_u+\theta^{(m)}_m,
\een
where each $\theta^{(k)}_m$ is a power series in $u$ with coefficients in
the algebra $\U(\wh\g^{\tss +})$. In particular,
\beql{trlm}
\theta^{(m)}_m=(-1)^m\tss\tr\ts\Lc^+(u)^m+\text{\ lower degree terms}.
\eeq
This follows from the expansion \eqref{tyexp}
and the relation
\ben
\tr^{}_{1,\dots,s}\ts P_{s-1\ts s}\dots P_{1\ts 2}\ts\Lc_1\dots \Lc_s=\tr\ts \Lc^s
=\tr\ts \big(2\tss u\tss \di_u-\Lc^+(u)\big)^s.
\een

Our main result is the following theorem which we will prove in the next section.

\bth\label{thm:commut}
The coefficients of all power series $\theta^{(k)}_m$ generate
a commutative subalgebra of $\U(\wh\g^{\tss +})$.
\eth

By Theorem~\ref{thm:commut} the commuting family \eqref{trlm} quantizes the well-known
Hamiltonians $\tr\ts L(u)^m$ of the classical trigonometric Gaudin model;
see \cite{bv:hs}, \cite{s:qi}. Note also that
the above expressions for $\theta_m$ with $m=1,2,3$ show that the commutative subalgebra
provided by Theorem~\ref{thm:commut} contains the coefficients
of the power series $\tr\ts \Lc^+(u)^2$ (see \eqref{quadham}),
as well as the coefficients of the power series
\ben
\tr\ts\Lc^+(u)^3-2\tss u\ts\tr\ts\Lc^+(u)\Lc^+(u)'
+\sum_{i,j=1}^N \sign(j-i)\Lc^+_{ij}(u)\Lc^+_{ji}(u).
\een

\section{Proof of Theorem~\ref{thm:commut}}
\label{sec:proof}

We start by recalling the {\em Bethe subalgebra} of the quantum affine algebra $\U_q(\wh\gl_N)$
over $\CC(q)$
associated with the $XXZ$ model; see e.g. \cite{rtv:tw} for a review.
This is a commutative subalgebra which lies within the $q$-{\em Yangian} $\Y_q(\gl_N)\subset \U_q(\wh\gl_N)$.
The algebra $\Y_q(\gl_N)$
is generated by elements
\ben
l^+_{ij}[-r],\qquad 1\leqslant i,j\leqslant N,\qquad r=0,1,\dots,
\een
with the conditions that $l^+_{ij}[0]=0$ for $i>j$ and the elements $l^+_{ii}[0]$ are invertible,
subject to the defining relations
\beql{rll}
R(u/v)L_1^{+}(u)L_2^{+}(v)=L_2^{+}(v)L_1^{+}(u)R(u/v).
\eeq
Here we use the matrix $L^{+}(u)=\big[\tss l^{+}_{ij}(u)\big]$,
whose entries are formal power series in $u$,
\ben
l^{+}_{ij}(u)=\sum_{r=0}^{\infty}l^{+}_{ij}[-r]\tss u^r
\een
and regard it as the element
\ben
L^{+}(u)=\sum_{i,j=1}^n e_{ij}\ot l^{+}_{ij}(u)\in\End\CC^N\ot\Y_q(\gl_N)[[u]].
\een
By a standard notation, subscripts are used to indicate copies of the matrix in the
tensor product algebra
\ben
\End\CC^N\ot\End\CC^N\ot\Y_q(\gl_N)[[u]]
\een
so that $L_2^{+}(v)=I\ot L^{+}(v)$ etc., where $I$ is the identity matrix.
The $R$-{\em matrix} is given by
\begin{align}
R(x)=\sum_{i}e_{ii}\ot e_{ii}
{}&+\frac{1-x}{q-\qin x}\ts\sum_{i\ne j}e_{ii}\ot e_{jj}
\non\\[0.4em]
{}+{}&\frac{(q-\qin)\ts x}{q-\qin x}\ts\sum_{i> j}e_{ij}\ot
e_{ji}+ \frac{q-\qin}{q-\qin x}\ts\sum_{i< j}e_{ij}\ot e_{ji}.
\label{Rx}
\end{align}

Consider the $q$-{\em permutation}
$P^{\tss q}\in\End(\CC^N\ot\CC^N)\cong\End\CC^N\ot\End\CC^N$
defined by
\beql{qperm}
P^{\tss q}=\sum_{i}e_{ii}\ot e_{ii}+ q\tss\sum_{i> j}e_{ij}\ot
e_{ji}+ \qin\sum_{i< j}e_{ij}\ot e_{ji}.
\eeq
The symmetric group $\Sym_k$ acts on the tensor product space $(\CC^N)^{\ot\tss k}$
by $s_a\mapsto P^{\tss q}_{s_a}:=
P^{\tss q}_{a\ts a+1}$ for $a=1,\dots,k-1$,
where $s_a$ denotes the transposition $(a,a+1)$. The operator
$P^{\tss q}_{a\ts a+1}$ acts as $P^{\tss q}$
in the tensor product of the $a$-th and $(a+1)$-th copies of $\CC^N$
and acts as the identity operator in
the remaining copies.
If $\si=s_{a_1}\dots s_{a_l}$ is a reduced decomposition
of an element $\si\in \Sym_k$ then we set
$P^{\tss q}_{\si}=P^{\tss q}_{s_{a_1}}\dots P^{\tss q}_{s_{a_l}}$.
We denote by $A^{(k)}$ the image of the normalized antisymmetrizer
associated with the $q$-permutations:
\beql{antisym}
A^{(k)}=\frac{1}{k!}\ts\sum_{\si\in\Sym_k}\sgn\tss\si\cdot P^{\tss q}_{\si}.
\eeq
For each $k=1,\dots,N$ consider the power series in $u$ defined by
\beql{lkz}
\tr^{}_{1,\dots,k}\ts A^{(k)}L^+_1(u)\dots L^+_k(uq^{-2k+2})
\eeq
with the trace taken over all $k$ copies of $\End\CC^N$ in the tensor product algebra
\beql{mtpy}
\underbrace{\End\CC^N\ot \dots\ot\End\CC^N}_{k}\ot\Y_q(\gl_N)[[u]].
\eeq
It is well-known that the coefficients of all power series \eqref{lkz}
generate a commutative subalgebra $\Bc_N$ of $\Y_q(\gl_N)$. Another family of generators of
this subalgebra can be obtained from the {\em Newton identities};
see \cite[Theorem~6.6]{cfrs:ap}. Adapting to our settings, we find that
the coefficients of all power series
\beql{plkz}
\tr^{}_{1,\dots,k}\ts P^{\tss q}_{(k,k-1,\dots,1)}L^+_1(u)\dots L^+_k(uq^{-2k+2}),\qquad k=1,2,\dots
\eeq
belong to $\Bc_N$.

Since the $q$-Yangian $\Y_q(\gl_N)$ is a deformation of the universal enveloping algebra
$\U(\wh\g^{\tss +})$,
the classical limit $q\to 1$ takes $\Bc_N$
to a commutative subalgebra of $\U(\wh\g^{\tss +})$. To get its generators
in an explicit form we will use the power series \eqref{plkz} and apply
an argument similar to the one used in \cite[Theorem~3.8]{km:cq}.
We will use both the permutation $P$ given in \eqref{permut}
and the $q$-permutation $P^{\tss q}$ defined in \eqref{qperm}.
Introduce the operator $\de$ which interacts with power series in $u$ by the rule
$\de\tss g(u)=g(uq^{-2})\tss\de$. Adjoining this element to the algebra $\Y_q(\gl_N)[[u]]$,
set $M=L^+(u)\tss \de$.
For each $m\geqslant 1$ consider the expression
\beql{temp_a}
\Mc_m=\frac{1}{(q-1)^m}\ts\big(1-(M_m)^{\rightarrow}  \big)
\Big(P^{}_{m-1\,m}-P^{\tss q}_{m-1\,m}(M_{m-1})^{\rightarrow}\Big)
\dots\Big(P^{}_{12}-P^{\tss q}_{12}(M_1)^{\rightarrow} \Big)\ts 1,
\eeq
where the arrow in the superscript indicates that the corresponding term
is understood as the operator of right multiplication:
\beql{apptx}
\Big(P^{}_{a\ts a+1}-P^{\tss q}_{a\ts a+1}\big(M_a \big)^{\rightarrow} \Big)\ts X:
=P^{}_{a\ts a+1}\ts X-P^{\tss q}_{a\ts a+1}\ts X\ts M_a.
\eeq
The operators in \eqref{temp_a} are meant to be applied consecutively from right to left. 
By taking the trace over all $m$ copies of $\End\CC^N$ in \eqref{temp_a}, we get
a polynomial in $\de$,
\ben
\tr^{}_{1,\dots,m}\ts\Mc_m\in \Y_q(\gl_N)[[u]][\de],
\een
whose coefficients are power series in $u$.

\ble\label{lem:mm}
All coefficients of the polynomial $\tr^{}_{1,\dots,m}\ts\Mc_m$
belong to $\Bc_N[[u]]$.
\ele

\bpf
Expand the product in \eqref{temp_a} to get the expression
\ben
\Mc_m=\frac{1}{(q-1)^m}\ts\sum_{k=0}^m \sum_{1\leqslant a_1<\dots< a_k\leqslant m }\ts(-1)^k\ts
\Pi_{a_1,\dots,a_k}\ts M_{a_1}\dots M_{a_k},
\een
where
\ben
\Pi_{a_1,\dots,a_k}=P^{}_{(m,m-1,\dots,a_k+1)}P^{\tss q}_{a_k\ts a_k+1}
P^{}_{(a_k,\dots,a_{k-1}+1)}P^{\tss q}_{a_{k-1}\ts a_{k-1}+1}\dots
P^{}_{(a_2,\dots,a_1+1)}P^{\tss q}_{a_1\ts a_1+1}P^{}_{(a_1,\dots,1)}.
\een
As shown in \cite[Theorem~3.8]{km:cq}, for the partial trace
we have
\beql{trpi}
\tr^{}_{\{1,\dots,m\}\setminus\{a_1,\dots, a_k\}}\ts \Pi_{a_1,\dots,a_k}
=P^{\tss q}_{a_{k-1}\ts a_{k}}P^{\tss q}_{a_{k-2}\ts a_{k-1}}\dots P^{\tss q}_{a_1\ts a_2}.
\eeq
Hence we obtain
\ben
\bal
\tr_{1,\dots,m}\tss\Mc_m&=
\frac{1}{(q-1)^m}\ts\sum_{k=0}^m \sum_{1\leqslant a_1<\dots< a_k\leqslant m }
\ts(-1)^k\ts\tr^{}_{a_1,\dots,a_k}\ts
P^{\tss q}_{a_{k-1}\ts a_{k}}
\dots P^{\tss q}_{a_1\ts a_2}\ts M_{a_1}\dots M_{a_k}\\[0.4em]
{}&=\frac{1}{(q-1)^m}\ts\sum_{k=0}^m \ts(-1)^k\ts \binom{m}{k}\tss \tr_{1,\dots,k}\ts
P^{\tss q}_{k-1\ts k}\dots P^{\tss q}_{1\ts 2}\ts M_1\dots M_k.
\eal
\een
Since $P^{\tss q}_{k-1\ts k}P^{\tss q}_{k-2\ts k-1}\dots P^{\tss q}_{1\ts 2}=P^{\tss q}_{(k,k-1,\dots,1)}$ and
\ben
M_1\dots M_k=
L^+_1(u)\dots L^+_k(uq^{-2k+2})\tss \de^k,
\een
the claim follows.
\epf

Lemma~\ref{lem:mm} provides a family of elements of the commutative subalgebra
$\Bc_N$ of $\Y_q(\gl_N)$. As a next step, we will calculate the classical limits $q\to 1$
of these elements. They will form a commuting family of elements of the algebra
$\U(\wh\g^{\tss +})$. To this end, we will use expansions into power series in $q-1$.
Write
\ben
\de=1-2\tss(q-1)\tss u\tss \di_u+\dots.
\een
We have
\ben
L^+(u)=1+(q-1)\tss\Lc^+(u)+\dots \Fand 1-M=(q-1)\tss\Lc+\dots
\een
with $\Lc=2\tss u\tss \di_u-\Lc^+(u)$.
Furthermore,
\ben
P-P^{\tss q}=(q-1)\tss\Tc+\dots
\een
with $\Tc$ defined in \eqref{permut}.
For the expression \eqref{apptx} we then get
\ben
\Big(P^{}_{a\ts a+1}-P^{\tss q}_{a\ts a+1}\big(M_a \big)^{\rightarrow} \Big)\ts X
=(q-1)\Big(\Tc_{a\ts a+1}+P^{}_{a\ts a+1}\big(\Lc_a \big)^{\rightarrow}\Big)\ts X+\dots.
\een

Thus we arrive at the next lemma.

\ble\label{lem:cllim}
The classical limit of the polynomial $\Mc_m$
is the differential operator
\beql{cllim}
\overline{\Mc}_m=(\Lc_m)^{\rightarrow}
\Big(\Tc^{}_{m-1\,m}+P_{m-1\,m}(\Lc_{m-1})^{\rightarrow}\Big)
\dots\Big(\Tc^{}_{12}+P_{12}(\Lc_1)^{\rightarrow} \Big)\ts 1,
\eeq
where we use the arrow notation as in \eqref{apptx}.
\qed
\ele

The expanded form of \eqref{cllim} is given by
\ben
\overline{\Mc}_m=\sum_{k=1}^m \sum_{1\leqslant a_1<\dots< a_k= m }\ts
\Ga_{a_1,\dots,a_k}\ts \Lc_{a_1}\dots \Lc_{a_k},
\een
where
\ben
\Ga_{a_1,\dots,a_k}=\Tc^{}_{(m,m-1,\dots,a_k+1)}P_{a_k\ts a_k+1}
\Tc^{}_{(a_k,\dots,a_{k-1}+1)}P_{a_{k-1}\ts a_{k-1}+1}\dots
\Tc^{}_{(a_2,\dots,a_1+1)}P_{a_1\ts a_1+1}\Tc^{}_{(a_1,\dots,1)}
\een
and we set
\ben
\Tc^{}_{(c_k,\dots,c_1)}=\Tc^{}_{c_{k-1}\ts c_k}\dots \Tc^{}_{c_1\ts c_2}
\een
and $\Tc^{}_{(a_i,\dots,a_{i-1}+1)}=1$ if $a_i=a_{i-1}+1$.

\ble\label{lem:trrma}
For any $k\geqslant 3$ we have
\ben
\tr^{}_{2,\dots,k-1}\tss\Tc^{}_{(k,k-1,\dots,1)}
=\begin{cases}
\Tc_{1\tss k}\qquad&\text{if}\quad $k$ \quad \text{is even},\\[0.4em]
\overline{\Tc}_{1\tss k}\qquad&\text{if}\quad $k$ \quad \text{is odd}.
\end{cases}
\een
\ele

\bpf
This follows by a straightforward induction argument.
\epf

As a next step, we use
Lemma~\ref{lem:trrma} and the identities
\ben
\tr^{}_1\Tc_{1\tss 2}=\tr^{}_1\overline{\Tc}_{1\tss 2}=0,
\een
to calculate the trace $\tr^{}_{1,\dots,m}\tss\overline{\Mc}_m$.
We also have
\ben
\tr^{}_2\tss \Tc_{2\tss 3}\tss P_{1\tss 2}=\tr^{}_2\tss P_{1\tss 2}\tss\Tc_{1\tss 3}=\Tc_{1\tss 3}
\een
and the same relation holds for $\Tc$ replaced by $\overline{\Tc}$.
We thus obtain
\ben
\tr^{}_{1,\dots,m}\tss\overline{\Mc}_m=
\sum_{k=1}^m \sum_{1=a_1<\dots< a_k= m }\ts
\tr^{}_{a_1,\dots,a_k}\tss\Tc^{}_{[a_{k-1}\ts a_k]}\dots \Tc^{}_{[a_1\ts a_2]}
\ts \Lc_{a_1}\dots \Lc_{a_k},
\een
where
\ben
\Tc^{}_{[a\ts b]}
=\begin{cases}
\Tc^{}_{a\ts b}\qquad&\text{if}\quad b-a\geqslant 3 \quad \text{is odd},\\[0.4em]
\overline{\Tc}_{a\ts b}\qquad&\text{if}\quad b-a\geqslant 2 \quad \text{is even},\\[0.4em]
P^{}_{a\ts b}\qquad&\text{if}\quad b-a=1.
\end{cases}
\een
We may thus conclude that the coefficients $\theta_m$ defined in
\eqref{thetham} are found by
\ben
\theta_m=\tr^{}_{1,\dots,m}\tss\overline{\Mc}_m.
\een
The coefficients of these differential operators pairwise commute
which completes the proof of Theorem~\ref{thm:commut}.

\paragraph{Shifted commutative subalgebra.}
To make a connection with invariants of the vacuum module over the Lie algebra
$\wh\g^{\tss +}$,
we also consider a modified version of the commutative subalgebra
of $\U(\wh\g^{\tss +})$ provided by Theorem~\ref{thm:commut}.
It is well-known that the relations \eqref{rll} remain valid after
the replacement $L^+(u)\mapsto L^+(u)\tss D$, where $D$ is the diagonal matrix
\beql{d}
D=\diag\big[q^{N-1}, q^{N-3},\dots, q^{-N+1}\big].
\eeq
Define differential operators $\vt_m\in \U(\wh\g^{\tss +})[[u,\di_u]]$ by means of the generating
function
\beql{vthetham}
\sum_{m=1}^{\infty} \vt_m\tss y^m=\sum_{s=1}^{\infty}\tss y^s\ts
\tr^{}_{1,\dots,s}\ts T_{s-1\ts s}(y)\dots T_{1\ts 2}(y)\ts\overline{\Lc}_1\dots \overline{\Lc}_s,
\eeq
where $\overline{\Lc}=2\tss u\tss \di_u-\rho-\Lc^+(u)$ and $\rho$ is the diagonal matrix
\beql{rho}
\rho=\diag\big[N-1, N-3,\dots, -N+1\big].
\eeq
The differential operator $\vt_m$ takes the form
\beql{vtm}
\vt_m=\vt^{(0)}_m\tss \di_u^{\tss m}+\dots+\vt^{(m-1)}_m\tss\di_u+\vt^{(m)}_m,
\eeq
where each $\vt^{(k)}_m$ is a power series in $u$ with coefficients in
the algebra $\U(\wh\g^{\tss +})$.
Repeating the arguments of this section for the matrix $M=L^+(u)\tss D\tss\de$,
we come to the following.

\bco\label{cor:vcommut}
The coefficients of the power series $\vt^{(k)}_m$ generate
a commutative subalgebra of $\U(\wh\g^{\tss +})$.
\eco

\bpf
The only additional observation is the power series expansion for the new matrix $M$
given by
\ben
M=1-(q-1)\tss \big(2\tss u\tss \di_u-\rho-\Lc^+(u)\big)+\dots
\een
implied by the expansion $D=1+(q-1)\tss \rho+\dots$.
\epf

\section{Invariants of the vacuum module}
\label{sec:invm}

Now we consider the full quantum affine algebra $\U_q(\wh\gl_N)$ in its
$RLL$ presentation; see \cite{fri:qa}, \cite{rs:ce}.
We will need
the normalized $R$-matrix
\beql{rf}
\overline R(x)=f(x)\tss R(x),
\eeq
where $R(x)$ is defined in \eqref{Rx} and
\ben
f(x)=1+\sum_{k=1}^{\infty}f_k(q)\tss x^k
\een
is a formal power series in $x$ whose coefficients $f_k(q)$ are
rational functions in $q$ uniquely determined
by the relation
\beql{fx}
f(xq^{2N})=f(x)\ts\frac{(1-xq^2)\tss(1-xq^{2N-2})}{(1-x)\tss(1-xq^{2N})}.
\eeq

The {\em quantum affine algebra $\U_q(\wh\gl_N)$}
is generated by elements
\ben
l^+_{ij}[-r],\qquad l^-_{ij}[r]\qquad\text{with}\quad 1\leqslant i,j\leqslant N,\qquad r=0,1,\dots,
\een
and the invertible central element $q^c$,
subject to the defining relations
\begin{align}
l^+_{ji}[0]&=l^-_{ij}[0]=0\qquad&&\text{for}\qquad 1\leqslant i<j\leqslant N,
\label{llt}\\
l^+_{ii}[0]\ts l^-_{ii}[0]&=l^-_{ii}[0]\ts l^+_{ii}[0]=1\qquad&&\text{for}\qquad i=1,\dots,N,
\label{lplm}
\end{align}
and
\begin{align}
R(u/v)L_1^{\pm}(u)L_2^{\pm}(v)&=L_2^{\pm}(v)L_1^{\pm}(u)R(u/v),
\label{RLL}\\[0.2em]
\overline R(uq^{-c}/v)L_1^{+}(u)L_2^{-}(v)&=L_2^{-}(v)L_1^{+}(u)\overline R(uq^{c}/v).
\label{RLLpm}
\end{align}
In the last two relations we consider the matrices $L^{\pm}(u)=\big[\tss l^{\pm}_{ij}(u)\big]$,
whose entries are formal power series in $u$ and $u^{-1}$,
\beql{serlpm}
l^{+}_{ij}(u)=\sum_{r=0}^{\infty}l^{+}_{ij}[-r]\tss u^r,\qquad
l^{-}_{ij}(u)=\sum_{r=0}^{\infty}l^{-}_{ij}[r]\tss u^{-r}.
\eeq
The $q$-Yangian $\Y_q(\gl_N)$ can be identified with the subalgebra of $\U_q(\wh\gl_N)$
generated by the coefficients of the series $l^{+}_{ij}(u)$ with $1\leqslant i,j\leqslant N$.

The {\em vacuum module at the critical level} $c=-N$ over $\U_q(\wh\gl_N)$
is the universal module $V_q(\gl_N)$ generated by a nonzero vector $\vac$ subject
to the conditions
\ben
L^-(u)\tss\vac=I\tss\vac,\qquad q^c\tss\vac=q^{-N}\tss\vac,
\een
where $I$ denotes the identity matrix.
As a vector space,
$V_q(\gl_N)$ can be identified with the subalgebra
$\Y_q(\gl_N)$ of $\U_q(\wh\gl_N)$ generated by
the coefficients of all series $l^+_{ij}(u)$
subject to the additional relations $l^+_{ii}[0]=1$.
The subspace of invariants of $V_q(\gl_N)$ is defined by
\ben
\z_q(\wh\gl_N)=\{v\in V_q(\gl_N)\ |\ L^-(u)\tss v=I v\}.
\een
One can regard $\z_q(\wh\gl_N)$ as a subspace of $\Y_q(\gl_N)$.
This subspace is closed under the multiplication in
the quantum affine algebra and it can be identified
with a subalgebra of $\Y_q(\gl_N)$. By \cite[Corollary~3.3]{fjmr:hs},
for $k=1,\dots,N$
all coefficients of the series
\beql{barell}
\tr^{}_{1,\dots,k}\ts A^{(k)}\tss L^+_1(z)\dots L^+_k(zq^{-2k+2})
D_1\dots D_k\tss\vac
\eeq
belong to the algebra of invariants $\z_q(\wh\gl_N)$.
Moreover, the coefficients of all these series pairwise commute.
As with the series \eqref{lkz}, applying the Newton identities
of \cite[Theorem~6.6]{cfrs:ap}, we find that
the coefficients of all power series
\beql{iplkz}
\tr^{}_{1,\dots,k}\ts P^{\tss q}_{(k,k-1,\dots,1)}L^+_1(u)\dots L^+_k(uq^{-2k+2})
D_1\dots D_k\tss\vac,
\eeq
belong to $\z_q(\wh\gl_N)$ for all $k\geqslant  1$.

Under the limit $q\to 1$ the algebra $\U_q(\wh\gl_N)$
turns into the universal enveloping
algebra $\U(\wh\gl_N)$. To make this statement more precise, consider
the presentation of the affine Lie algebra $\wh\gl_N$ associated with the classical
$r$-matrix \eqref{rx}. Recall that
the affine Kac--Moody algebra $\wh\gl_N=\gl_N[t,t^{-1}]\oplus\CC K$
has the commutation relations
\beql{commrel}
\big[E_{ij}[r],E_{kl}[s\tss]\tss\big]
=\de_{kj}\ts E_{i\tss l}[r+s\tss]
-\de_{i\tss l}\ts E_{kj}[r+s\tss]
+r\tss\de_{r,-s}\ts K\Big(\de_{kj}\tss\de_{i\tss l}
-\frac{\de_{ij}\tss\de_{kl}}{N}\Big),
\eeq
and the element $K$ is central, where we set $E_{ij}[r]=E_{ij}\tss t^r$
for all $r\in\ZZ$. In addition to the matrix
$\Lc^+(u)$ with the entries
\eqref{lplus} introduce the matrix
$\Lc^-(u)=[\Lc^-_{ij}(u)]$ with
\beql{lminus}
\Lc^-_{ij}(u)=\sum_{n=0}^{\infty} \Lc^-_{ij}[n]\tss u^{-n},
\eeq
where $\Lc^-_{ij}[n]=2\tss E_{ij}[n]$ for $n\geqslant 1$ and
$\Lc^-_{ij}[0]=\big(1+\sign(i-j)\big) E_{ij}[0]$.
The defining relations of the algebra $\U(\wh\gl_N)$ can be written in the form
\begin{align}
\label{lpmpm}
\big[\Lc^{\pm}_1(u),\Lc^{\pm}_2(v)\big]&=\big[\Lc^{\pm}_1(u)+\Lc^{\pm}_2(v),r_{12}(u/v)\big],\\[0.4em]
\big[\Lc^{+}_1(u),\Lc^{-}_2(v)\big]&=\big[\Lc^{+}_1(u)+\Lc^{-}_2(v),r_{12}(u/v)\big]
+\frac{4 uv}{(u-v)^2}\Big(P_{12}-\frac{1}{N}\Big)\tss K,
\label{lmlp}
\end{align}
where $r(x)$ is defined in \eqref{rx} and we write $1$ for the tensor product
of the identity matrices $I\ot I$; cf. \cite{bbc:fp}.
We have the following
well-known
property.

\bpr\label{prop:cllim}
The defining relations of $\U(\wh\gl_N)$ are recovered from those of
$\U_q(\wh\gl_N)$ by the expansions into power series in $q-1$,
\ben
L^{\pm}(u)=I+(q-1)\tss \Lc^{\pm}(u)+\dots
\een
and setting $c\mapsto K$.
\epr

\bpf
We will only demonstrate how the relation \eqref{lmlp} is obtained
from \eqref{RLLpm}, which should explain the role of the normalized $R$-matrix \eqref{rf}.
Relations \eqref{lpmpm} are verified in the same way with simpler calculations.
Expanding into power series in $q-1$ and identifying $I\ot I$ with $1$ we get
\beql{Rxexp}
R(x)=1+(q-1)\tss \Big(r(x)-\frac{1+x}{1-x}\ts 1\Big)+\dots
\eeq
and
\beql{fxexp}
f(x)=1+2\tss (q-1)\tss\frac{(N-1)\tss x}{N\tss(1-x)}+\dots,
\eeq
where the second expansion is implied e.g. by the calculations in
\cite[Sec.~2]{km:cq}. Now apply \eqref{RLLpm} to get
\begin{multline}
\overline R(uq^{-c}/v)\big(L_1^{+}(u)-1\big)\big(L_2^{-}(v)-1\big)-
\big(L_2^{-}(v)-1\big)\big(L_1^{+}(u)-1\big)\tss\overline R(uq^{c}/v)\\[0.5em]
{}=\big(L_1^{+}(u)-1+L_2^{-}(v)-1\big)\tss\overline R(uq^{c}/v)
-\overline R(uq^{-c}/v)\big(L_1^{+}(u)-1+L_2^{-}(v)-1\big)\\[0.5em]
{}+\overline R(uq^{c}/v)-\overline R(uq^{-c}/v).
\non
\end{multline}
Dividing
both sides by $(q-1)^2$ and taking the limit $q\to 1$ we get
\ben
\big[\Lc^{+}_1(u),\Lc^{-}_2(v)\big]=\big[\Lc^{+}_1(u)+\Lc^{-}_2(v),r_{12}(u/v)\big]
+\frac{\overline R(uq^{c}/v)-\overline R(uq^{-c}/v)}{(q-1)^2}\Big|_{q=1}.
\een
Using \eqref{Rxexp} and \eqref{fxexp} we find that
\ben
\frac{\overline R(xq^{c})-\overline R(xq^{-c})}{(q-1)^2}\Big|_{q=1}
=\frac{4\tss c\tss x}{(1-x)^2}\Big(P-\frac{1}{N}\Big)
\een
thus completing the proof.
\epf

The ({\em trigonometric}) {\em vacuum module at the critical level} over the affine Lie algebra $\wh\gl_N$
is the universal module $V^{}_{\text{\rm tr}}(\gl_N)$ generated by a nonzero vector $\vac$ subject
to the conditions $K\tss\vac=-N\tss\vac$ and
\beql{vacdef}
E_{ij}[n]\tss\vac=0\quad \text{for all $i,j$}\fand n\geqslant 1,\Fand
E_{ij}[0]\tss\vac=0\quad \text{for $i\geqslant j$}.
\eeq
By \eqref{lminus} these conditions can be written in a matrix form as
$\Lc^-(u)\tss\vac=0$.
The {\em subspace of invariants of} $V^{}_{\text{\rm tr}}(\gl_N)$ is defined by
\beql{invdef}
\z^{}_{\text{\rm tr}}(\wh\gl_N)=\{v\in V^{}_{\text{\rm tr}}(\gl_N)\ |\ \Lc^-(u)\tss v=0\}.
\eeq
By the Poincar\'e--Birkhoff--Witt theorem,
the vacuum module is isomorphic to the universal enveloping algebra
$\U(\wh\g^{\tss +})$, as a vector space, so that
we can regard $\z^{}_{\text{\rm tr}}(\wh\gl_N)$ as a subalgebra of $\U(\wh\g^{\tss +})$.

Note that the above definitions are quite analogous to those of the standard vacuum module
over $\wh\gl_N$, where the conditions \eqref{vacdef} are replaced by
\ben
E_{ij}[n]\tss\vac=0\quad \text{for all $i,j$}\fand n\geqslant 0.
\een
The corresponding vacuum module has a vertex algebra structure and
its subspace of invariants defined by analogy with \eqref{invdef}
coincides with the center of this vertex algebra. The center is a commutative associative
algebra whose structure was described by a theorem of Feigin and Frenkel~\cite{ff:ak};
see also \cite{f:lc} for a detailed general proof and \cite{m:so}
for an explicit approach in the case of classical Lie algebras.
Our next result can be regarded as a trigonometric analogue
of the Sugawara operators in type $A$. Recall
the power series $\vt^{(k)}_m\in \U(\wh\g^{\tss +})[[u]]$ defined in
\eqref{vthetham} and \eqref{vtm}.

\bth\label{thm:comminv}
The coefficients of the power series $\vt^{(k)}_m$
belong to $\z^{}_{\text{\rm tr}}(\wh\gl_N)$.
\eth

\bpf
As we pointed out above,
the coefficients of all power series
\eqref{iplkz} are invariants of the vacuum module $V_q(\gl_N)$ over
the quantum affine algebra $\U_q(\wh\gl_N)$. Therefore, the claim is derived from
Proposition~\ref{prop:cllim}
by taking the limit $q\to 1$ as in Section~\ref{sec:proof}.
\epf

\end{document}